\documentclass[11pt]{amsart}
\usepackage[T1]{fontenc}
\usepackage{lmodern}
\usepackage[a4paper,margin=2cm]{geometry}
\pdfpagewidth=\paperwidth \pdfpageheight=\paperheight
\usepackage{amsfonts,amssymb,amsthm,amsmath,eucal}
\usepackage{ltablex}
\usepackage{caption, booktabs}
\usepackage{longtable}
\usepackage{hyperref}
\theoremstyle{plain}
\newtheorem{thm}{Theorem}[section]
\newtheorem{theorem}[thm]{Theorem}
\newtheorem*{theoremA}{Theorem A}
\newtheorem*{theoremB}{Theorem B}
\newtheorem*{theoremC}{Theorem C}
\newtheorem{lemma}[thm]{Lemma}
\newtheorem{proposition}[thm]{Proposition}
\newtheorem{corollary}[thm]{Corollary}

\theoremstyle{definition}
\newtheorem{remark}[thm]{Remark}

\newcommand\be{\begin{eqnarray*}}
\newcommand\ee{\end{eqnarray*}}

\newcommand\C{\mathbb C}
\DeclareMathOperator{\GL}{GL}
\newcommand\call{\mathcal L}
\newcommand\cali{\mathcal I}

\newcommand\Z{\mathbb Z}

\newcommand\cH{\mathcal H}

\newcommand{\Loc}{\operatorname{Loc}}
\newcommand\rationalto{\dashrightarrow}

\renewcommand\P{\mathbb P}
\renewcommand\O{\mathcal O}
\newcommand\newop[2]{\def#1{\mathop{\rm #2}\nolimits}}
\newop\SL{SL}
\newop\mult{mult}
\newop\alphahat{\widehat{\alpha}}

\begin{document}

\author[K. Hanumanthu]{Krishna Hanumanthu}
\address{Chennai Mathematical Institute, H1 SIPCOT IT Park, Siruseri, Kelambakkam 603103, India}
\email{krishna@cmi.ac.in}

\author[J. Huizenga]{Jack Huizenga}
\address{Department of Mathematics, The Pennsylvania State University, University Park, PA 16802, USA}
\email{huizenga@psu.edu}

\author[J. Kettinger]{Jake Kettinger}
\address{Department of Mathematics, Colorado State University, Fort Collins, CO 80523, USA}
\email{jake.kettinger@colostate.edu}

\author[Ł. Merta]{Łukasz Merta}
\address{Department of Mathematics, University of the National Education Commission, PL-30-084 Krakow, Poland}
\email{lukasz.merta.26@gmail.com}

\author[T. Szemberg]{Tomasz Szemberg}
\address{Department of Mathematics, University of the National Education Commission, PL-30-084 Krakow, Poland}
\email{tomasz.szemberg@gmail.com}

\author[J. Szpond]{Justyna Szpond}
\address{Department of Mathematics, University of the National Education Commission, PL-30-084 Krakow, Poland}
\email{szpond@gmail.com}

\title[Waldschmidt constants for exceptional rank-four complex reflection groups]{Waldschmidt constants and invariant geometry for exceptional rank-four complex reflection groups}
\date{\today}

\thispagestyle{empty}
\begin{abstract}
We compute the Waldschmidt constants of the point configurations dual to the reflection arrangements of the exceptional complex reflection groups $G_{29}$, $G_{31}$ and $G_{32}$: they equal $3$ for the forty points of $G_{29}$, $36/11$ for the sixty-point Klein configuration, and $10/3$ for the forty-point Penrose configuration.
Together with the known cases $W(F_4)$ and $W(H_4)$, this completes the computation for all exceptional irreducible complex reflection groups of rank four.
An invariant dodecic and its sixfold orbit product give extremal surfaces for the first two configurations, and we identify the six dodecics as pullbacks of the coordinate hyperplane sections of the Segre cubic under the Igusa--Segre duality.
The matching lower bounds follow from restriction to quadrics and distinguished planes, including a weighted B\'ezout estimate arising from two related Hesse arrangements.
As a consequence, we determine the asymptotic resurgences of the three ideals.
\end{abstract}
\maketitle

\section{Introduction}
Let $Z\subset \P^N$ be a finite set of points and let $I=I(Z)$ be its homogeneous ideal.  
The $m$-th symbolic power $I^{(m)}$ consists of the forms vanishing to order at least $m$ at every point of $Z$.
Writing $\alpha(J)$ for the least degree of a non-zero homogeneous form in a homogeneous ideal $J$, 
the \emph{Waldschmidt constant} of $Z$ is
$$
\alphahat(Z)=\lim_{m\to\infty}\frac{\alpha(I^{(m)})}{m}
=\inf_{m\geq 1}\frac{\alpha(I^{(m)})}{m}.
$$
Thus $\alphahat(Z)$ measures the asymptotic cost in degree of imposing uniformly increasing multiplicity along $Z$.  
Waldschmidt constants occur naturally in the study of symbolic powers and containment problems, and they also admit a birational interpretation in terms of the boundary of the effective cone of the blow-up of projective space at the points; 
see for example \cite{BocciHarbourne2010,DumnickiSzembergSzpond2024}.

Highly symmetric point configurations provide a particularly attractive testing ground for this invariant.
A finite complex reflection group $G\subset\mathrm{GL}(V)$ determines an arrangement $\mathcal A(G)$ of reflecting hyperplanes and, by ordinary projective duality, a finite configuration
$$
Z_G\subset\mathbb P(V^*).
$$
For the computations, we choose a $G$-invariant positive definite Hermitian form and work in orthonormal coordinates.
We represent the dual configuration by the normal directions to the reflecting hyperplanes in $\mathbb P(V)$ and use the same notation $Z_G$ for this realization.
The identification with the dual configuration is conjugate-linear; conjugating the coefficients of a form preserves its degree and its orders of vanishing at conjugate points, so the Waldschmidt constant is unchanged.
For a point $P=[v]\in\mathbb P(V)$, with $v$ a column vector, we write $H_P=P^\perp=\mathbb P(v^\perp)$, where $v^\perp=\{X\in V:\overline v^{\,T}X=0\}$.
More generally, if $L=\mathbb P(U)$ is a projective linear subspace, we put $L^\perp=\mathbb P(U^\perp)$.
Unitarity gives $H_{gP}=g(H_P)$ for $g\in G$.

Such configurations often carry substantially more geometry than is visible from the defining arrangement alone.
They have appeared in connection with symbolic powers, unexpected hypersurfaces, special projections, and \emph{geproci} sets.
The Klein configuration of $60$ points is a particularly rich example \cite{PokoraSzembergSzpond2024}.
From the invariant-theoretic point of view, reflection groups are especially accessible because the Shephard--Todd--Chevalley theorem makes the global invariant ring polynomial \cite{ShephardTodd1954,Chevalley1955}, while Steinberg's theorem gives the same structure for stabilizers of vectors \cite{Steinberg1964}.

The Waldschmidt constants of the configurations in $\P^3$ arising from the root systems $D_4$, $B_4$, $F_4$ and $H_4$ were computed in \cite{CalvoHuizengaSzemberg2026}.
Among the exceptional irreducible complex reflection groups of rank four, these results cover $G_{28}=W(F_4)$ and $G_{30}=W(H_4)$.
Here we treat the remaining groups $G_{29}$, $G_{31}$ and $G_{32}$, combining invariant surfaces with restriction arguments adapted to the incidence geometry of their configurations.
General background on finite complex reflection groups can be found in \cite{LehrerTaylor2009,Kane2001}.

Central to our approach is the inclusion $G_{29}\subset G_{31}$ of index six \cite{BonnafeSarti2021}.
Let $Z_{29}\subset\P^3$ denote the configuration of $40$ points associated with $G_{29}$.
We identify a $G_{29}$-invariant dodecic surface $S_{12}=V(B_{12})$ whose equation vanishes to order $4$ at every point of $Z_{29}$.
The product of defining equations of its six $G_{31}$-translates yields an extremal surface for the Klein configuration associated with $G_{31}$.
Thus a common dodecic and its orbit product govern the two computations.
We explain this construction geometrically through the classical relation between the Igusa quartic and the Segre cubic.

The matching lower bounds use different restriction geometries.
For $G_{29}$ we restrict to distinguished planes containing twelve configuration points, whereas for the Klein configuration we use the grids on its fundamental quadrics.
For $G_{32}$, restriction to polar planes produces two related Hesse arrangements whose incidences yield a weighted B\'ezout estimate.
Together with estimates for vanishing along lines, these restrictions give the sharp lower bounds.
In particular, for $Z_{29}$ the lower bound matches the ratio $12/4=3$ supplied by $S_{12}$, giving our first result.
\begin{theoremA}
For the $G_{29}$ configuration $Z_{29}\subset\P^3$ one has
$$
\alphahat(Z_{29})=3.
$$
\end{theoremA}
Now let $Z_{31}$ be the $60$-point Klein configuration associated with $G_{31}$ and studied in \cite{PokoraSzembergSzpond2024}.
The set $Z_{31}$ contains $Z_{29}$ as a subset.
We show that $B_{12}$ vanishes to order $3$ at every point of $Z_{31}\setminus Z_{29}$.
For representatives $g_1,\ldots,g_6$ of the six left cosets of $G_{29}$ in $G_{31}$, the forms $g_i\cdot B_{12}$ define the six translates $g_i(S_{12})$.
Their product
$$
B_{72}=\prod_{i=1}^6(g_i\cdot B_{12})
$$
is a $G_{31}$-invariant form of degree $72$ vanishing uniformly to order $22$ at the points of $Z_{31}$.
The reducible surface $S_{72}=V(B_{72})$ yields our second principal result.
\begin{theoremB}
For the Klein configuration $Z_{31}\subset\P^3$ one has
$$
\alphahat(Z_{31})=\frac{36}{11}.
$$
\end{theoremB}
The dodecic $S_{12}$ has an intrinsic geometric description through the classical duality between the Igusa quartic and the Segre cubic.
Let $H_{2,2}\subset \mathrm{SL}_4(\C)$ be the finite Heisenberg group of order $32$, a central extension of $(\Z/2\Z)^4$ by $\{\pm I\}$. 
We use the realization and coordinates described in \cite{PokoraSzembergSzpond2024}.
Then the Heisenberg-invariant quartics define the Maschke map
$$
\mu:\P^3\longrightarrow \mathcal I\subset \P^4,
$$
whose image is the Igusa quartic $\mathcal I$.  
The Gauss map $\gamma$ of $\mathcal I$ maps onto its projective dual, the Segre cubic $\mathcal S\subset (\P^4)^\vee$; see \cite{Hunt1996,Dolgachev2012}.
We show that $S_{12}$ is the pullback, under the composite map
$$
\P^3 \stackrel{\mu}{\longrightarrow} \mathcal I \stackrel{\gamma}{\dashrightarrow} \mathcal S,
$$
of one of six distinguished hyperplane sections of the Segre cubic. The other five sections pull back to the remaining $G_{31}$-translates of $S_{12}$. Thus the classical Igusa--Segre duality explains simultaneously the occurrence of the $G_{29}$-invariant dodecic and the sixfold orbit product which yields $S_{72}$.

Finally, we turn to the remaining exceptional rank-four group $G_{32}$.  Here the relevant set is the $40$-point Penrose configuration $Z_{32}\subset\P^3$, originating in the work of Zimba and Penrose on Bell nonlocality \cite{ZimbaPenrose1993}. 
From the algebraic-geometric point of view, $Z_{32}$ is also a highly special $(5,8)$-geproci set \cite{ChiantiniEtAl2022}.
Restricting the reflection arrangement to one of the forty reflecting hyperplanes of $G_{32}$ lowers the symmetry to $G_{25}$ and produces, in the plane, the Hesse configuration together with its Hermitian polar dual; the relevant reflection-arrangement geometry is described in \cite{CousinMoussard2018}.  
The incidences of these two planar arrangements yield a weighted B\'ezout estimate, which combines with vanishing along four-point lines to give the sharp lower bound.
\begin{theoremC}
Let $Z_{32}\subset\P^3$ be the $40$-point Penrose configuration. Then
$$
\alphahat(Z_{32})=\frac{10}{3}.
$$
\end{theoremC}
Together with the results for $G_{28}=W(F_4)$ and $G_{30}=W(H_4)$ in \cite{CalvoHuizengaSzemberg2026}, Theorems~A--C determine the Waldschmidt constants of the point configurations associated with all five exceptional irreducible complex reflection groups of rank four.
Table~\ref{tab:exceptional-waldschmidt} summarizes these values.

\begin{table}[ht]
\centering
\caption{Waldschmidt constants for exceptional rank-four complex reflection groups.}
\label{tab:exceptional-waldschmidt}
\begin{tabular}{lccc}
\toprule
Group & $|Z_G|$ & $\alphahat(Z_G)$ & Reference \\
\midrule
$G_{28}=W(F_4)$ & $24$ & $8/3$
  & \cite{CalvoHuizengaSzemberg2026} \\
$G_{29}$ & $40$ & $3$ & Theorem~A \\
$G_{30}=W(H_4)$ & $60$ & $18/5$
  & \cite{CalvoHuizengaSzemberg2026} \\
$G_{31}$ (Klein) & $60$ & $36/11$ & Theorem~B \\
$G_{32}$ (Penrose) & $40$ & $10/3$ & Theorem~C \\
\bottomrule
\end{tabular}
\end{table}
In Section~\ref{sec:resurgence} we deduce the asymptotic resurgences of the three ideals and bounds for their resurgences (Corollary~\ref{cor:resurgence}).

\section{Invariant surfaces and B\'ezout estimates}
\label{sec:invariant-surfaces}

The configurations studied in this paper arise as projective orbits of finite complex reflection groups. 
We first explain how this symmetry can be used to study surfaces with prescribed multiplicities at their points.

Let $G\subset\GL_4(\C)$ be a finite group acting on the ambient space of a point configuration. 
For configurations arising from reflection arrangements, we use the unitary reflection representation in which the configuration consists of normal directions, as fixed in the introduction.
Let $S=\C[x,y,z,w]$
be the homogeneous coordinate ring, with the induced action
$(g\cdot F)(v)=F(g^{-1}v)$, and let
$$
T=S^G=\bigoplus_{d\geq0}T_d
$$
be its ring of invariants. 
For a complex reflection group the Shephard--Todd--Chevalley theorem gives
$$
T=\C[f_{d_1},f_{d_2},f_{d_3},f_{d_4}],
\mbox{ with } \deg f_{d_i}=d_i,
$$
where the generators are homogeneous and algebraically independent
\cite{ShephardTodd1954,Chevalley1955}.

Let $p\in\P^3$ be a point. We denote by $Z=G\cdot p$ its orbit. We consider two homogeneous ideals in $S$: that of the point $p$ denoted by $\mathfrak p=I(p)$ and that of the orbit $Z$ denoted by $J=I(Z)$.
An invariant form has the same order of vanishing at every point of $Z$.
Consequently, for $F\in T_d$ and $m\geq1$ there is the equivalence
$$
F\in J^{(m)}
\quad\Longleftrightarrow\quad
F\in\mathfrak p^m.
$$
Thus uniform vanishing along the orbit can be checked at one point.
The following proposition shows that restricting to invariant forms
also preserves the optimal asymptotic ratio of degree to multiplicity.

\begin{proposition}[Invariant reduction]\label{prop:invariant-reduction}
The Waldschmidt constant of a single orbit can be computed from $G$-invariant forms
$$
\alphahat(Z)=
\inf\left\{\frac{d}{m}\;:\; 
0\neq F\in T_d\cap\mathfrak p^m,\ m\geq 1 \right\}.
$$
\end{proposition}

\begin{proof}
Every non-zero form in $T_d\cap\mathfrak p^m$ belongs to $J^{(m)}$, so the right-hand side is at least $\alphahat(Z)$. 
To see the reverse inequality, let
$0\neq F\in J^{(m)}_d$. 
Its norm
$N_G(F)=\prod_{g\in G}g\cdot F$
is a non-zero $G$-invariant form of degree $|G|d$. 
Since every factor belongs to $J^{(m)}$, the product belongs to $J^{(|G|m)}$. Thus $N_G(F)$ has the same or a smaller degree to multiplicity ratio as $F$.
\end{proof}
Proposition~\ref{prop:invariant-reduction} justifies our search for extremal equations among invariant forms when constructing upper bounds.
The matching lower bounds will be proved for arbitrary surfaces by restriction to distinguished planes and quadrics.

We will also use the following elementary consequence of B\'ezout's theorem.
It bounds the multiplicity along a line in terms of prescribed multiplicities at points on that line.

\begin{lemma}\label{lem:penrose-forced-line}
Let $C$ be an effective plane divisor with rational coefficients and degree $c$, and let $L$ be a line containing $s\geq2$ distinct points at which $C$ has multiplicities at least $\mu_1,\ldots,\mu_s$.
Then the coefficient of $L$ in $C$ is at least $(\mu_1+\cdots+\mu_s-c)/(s-1)$.
The same estimate holds for the order of vanishing of a surface along a line containing the prescribed points.
\end{lemma}
\begin{proof}
We write $C=qL+C'$, with $C'$ not containing $L$.
The degree of $C'$ is $c-q$, and its multiplicity at the $i$th point is at least $\mu_i-q$.
B\'ezout gives $c-q\geq\sum_i(\mu_i-q)$, which is the required estimate.

For a surface, restrict to a general plane containing the line.
The coefficient of the line in this plane section equals the order of vanishing of the surface along the line, so the same argument applies.
\end{proof}

\section{The Klein group and its \texorpdfstring{$G_{29}$}{G29} subgroup}

\label{sec:klein-g29}

We now apply the invariant-theoretic framework to the inclusion
$G_{29}\subset G_{31}$. We use homogeneous coordinates $[x:y:z:w]$ on
$\P^3$ and write $Z_{29}\subset Z_{31}$ for the corresponding
configurations of forty and sixty points introduced above. A bar denotes
the projective image of a linear group. We first describe their common
line and quadric geometry, and then identify the invariant dodecic through
the classical duality between the Igusa quartic and the Segre cubic.

\subsection{The index-six subgroup and the point configurations}

Let $G=G_{31}\subset\GL_4(\C)$ be the Klein reflection group in the Shephard--Todd notation. 
Its order is $46080$, its scalar centre is $\mu_4$, and its basic invariant degrees are $8,12,20,24$.
Thus $|\overline G|=11520$.
We use the coordinates and numbering of \cite[Section 3, p.~603]{PokoraSzembergSzpond2024}. The set denoted there by $Z_{60}$ is our $Z_{31}=\{P_1,\ldots,P_{60}\}$:

\begingroup\small
\setlength{\tabcolsep}{5pt}
\begin{longtable}{lll}
$P_{1}=[0:0:1:1]$ & $P_{2}=[0:0:1:i]$ & $P_{3}=[0:0:1:-1]$ \\
$P_{4}=[0:0:1:-i]$ & $P_{5}=[0:1:0:1]$ & $P_{6}=[0:1:0:i]$ \\
$P_{7}=[0:1:0:-1]$ & $P_{8}=[0:1:0:-i]$ & $P_{9}=[0:1:1:0]$ \\
$P_{10}=[0:1:i:0]$ & $P_{11}=[0:1:-1:0]$ & $P_{12}=[0:1:-i:0]$ \\
$P_{13}=[1:0:0:1]$ & $P_{14}=[1:0:0:i]$ & $P_{15}=[1:0:0:-1]$ \\
$P_{16}=[1:0:0:-i]$ & $P_{17}=[1:0:1:0]$ & $P_{18}=[1:0:i:0]$ \\
$P_{19}=[1:0:-1:0]$ & $P_{20}=[1:0:-i:0]$ & $P_{21}=[1:1:0:0]$ \\
$P_{22}=[1:i:0:0]$ & $P_{23}=[1:-1:0:0]$ & $P_{24}=[1:-i:0:0]$ \\
$P_{25}=[1:0:0:0]$ & $P_{26}=[0:1:0:0]$ & $P_{27}=[0:0:1:0]$ \\
$P_{28}=[0:0:0:1]$ & $P_{29}=[1:1:1:1]$ & $P_{30}=[1:1:1:-1]$ \\
$P_{31}=[1:1:-1:1]$ & $P_{32}=[1:1:-1:-1]$ & $P_{33}=[1:-1:1:1]$ \\
$P_{34}=[1:-1:1:-1]$ & $P_{35}=[1:-1:-1:1]$ & $P_{36}=[1:-1:-1:-1]$ \\
$P_{37}=[1:1:i:i]$ & $P_{38}=[1:1:i:-i]$ & $P_{39}=[1:1:-i:i]$ \\
$P_{40}=[1:1:-i:-i]$ & $P_{41}=[1:-1:i:i]$ & $P_{42}=[1:-1:i:-i]$ \\
$P_{43}=[1:-1:-i:i]$ & $P_{44}=[1:-1:-i:-i]$ & $P_{45}=[1:i:1:i]$ \\
$P_{46}=[1:i:1:-i]$ & $P_{47}=[1:-i:1:i]$ & $P_{48}=[1:-i:1:-i]$ \\
$P_{49}=[1:i:-1:i]$ & $P_{50}=[1:i:-1:-i]$ & $P_{51}=[1:-i:-1:i]$ \\
$P_{52}=[1:-i:-1:-i]$ & $P_{53}=[1:i:i:1]$ & $P_{54}=[1:i:-i:1]$ \\
$P_{55}=[1:-i:i:1]$ & $P_{56}=[1:-i:-i:1]$ & $P_{57}=[1:i:i:-1]$ \\
$P_{58}=[1:i:-i:-1]$ & $P_{59}=[1:-i:i:-1]$ & $P_{60}=[1:-i:-i:-1]$ \\
\end{longtable}
\endgroup
We fix the forty-point subconfiguration
\begin{equation}\label{eq:z29-coordinates}
Z_{29}=\{P_1,\ldots,P_{24}\}\cup\{P_j:j\in J_+\},
\end{equation}
with
$ J_+=\{29,32,34,35,38,39,41,44,46,47,49,52,54,55,57,60\}$.
Equivalently, $J_+$ selects the points among $P_{29},\ldots,P_{60}$ whose coordinate product is $1$. 
This condition refers to the normalized
representatives in the table, whose first coordinate is $1$.
The complementary set $R_{20}=Z_{31}\setminus Z_{29}$ consists of $P_{25},\ldots,P_{28}$ and the sixteen points with coordinate product $-1$.

Selecting product $-1$ instead gives a conjugate embedding of $G_{29}$: the change $w\mapsto-w$ interchanges the two forty-point configurations.
The comparison with the dodecic displayed by Bonnaf\'e will be made through the explicit polynomial identity \eqref{eq:bonnafe-comparison}.
We continue to use the product equal $1$ convention throughout this section.

For each point $P_j=[a_j:b_j:c_j:d_j]$ in the table above, let $p_j=(a_j,b_j,c_j,d_j)$ be its coordinate row vector.
We index the reflections by coefficient vectors and define
$$
r_j=I_4-2\frac{\overline{p_j}^{\,T}p_j}
                 {p_j\overline{p_j}^{\,T}}.
$$
Here $I_4$ is the $4\times4$ identity matrix, and overlines on vectors denote complex conjugation.
This reflection fixes the hyperplane
$$
H_{\overline{P_j}}=V(a_jx+b_jy+c_jz+d_jw)
$$
pointwise and acts as multiplication by $-1$ on the line spanned by $\overline{p_j}^{\,T}$.
The formula is independent of the chosen homogeneous representative of $P_j$.
Thus $r_j$ is the reflection whose Hermitian normal direction is $\overline{P_j}$.

Complex conjugation preserves each of $Z_{31}$, $Z_{29}$ and $R_{20}$.
Consequently, on each of these sets, the families $\{H_{\overline P}\}$ and $\{H_P\}$ coincide as unlabelled families of planes.
We write $r_j=r_{P_j}$ with the indexing just specified.
The reflections with indices in \eqref{eq:z29-coordinates} generate our chosen subgroup $K\subset G$ isomorphic to $G_{29}$.
More specifically,
\begin{equation}\label{eq:generators-of-K}
    K=\langle r_{23},r_{11},r_{22},r_{29}\rangle.
\end{equation}
The first three generators interchange $x,y$, interchange $y,z$, and send
$(x,y,z,w)$ to $(-iy,ix,z,w)$, 
respectively. 
The fourth is
$$
r_{29}=\frac12\begin{pmatrix}
1&-1&-1&-1\\
-1&1&-1&-1\\
-1&-1&1&-1\\
-1&-1&-1&1
\end{pmatrix}.
$$
This subgroup has order $7680$, centre $\mu_4$, and basic invariant degrees $4,8,12,20$. 
Its projective action has two orbits on $Z_{31}$, namely $Z_{29}$ and $R_{20}$, and
$$
[G:K]=[\overline G:\overline K]=6,\;\; |\overline K|=1920.
$$
The degree-four invariant for this embedding is
\begin{equation}\label{eq:g29-quartic}
 a_4=x^4+y^4+z^4+w^4-12xyzw.
\end{equation}
Direct substitution in $a_4$ shows that precisely the forty reflections specified in \eqref{eq:z29-coordinates} preserve it. 
Adjoining $r_{25}$ to the four generators displayed in \eqref{eq:generators-of-K} gives $G$.
Each of the forty corresponding planes contains nine points of $Z_{29}$, whereas each of the twenty remaining Klein planes contains twelve points of $Z_{29}$.
Thus the forty planes and forty points form a $(40_9)$ configuration.

\subsection{The thirty lines and the ten fundamental quadrics}

Let $\call$ be the set of thirty Klein lines. 
Six are the coordinate lines. 
The other twenty-four are obtained by partitioning the coordinates into two pairs and imposing
$$
x_a=\epsilon x_b,\;\; x_c=\eta x_d,\;\;
(\epsilon,\eta)\in\{1,-1\}^2\ \cup\ \{i,-i\}^2.
$$
There are three such partitions and eight choices for each. The lines meet three at a time at the points of $Z_{31}$, and each line contains six such points \cite[Sections 2--3]{PokoraSzembergSzpond2024}. 
In the subdivision above, each line contains four points of $Z_{29}$ and two points of $R_{20}$. 
In particular, the same thirty lines support both point configurations.

We keep the following numbering for the fundamental quadrics:
$$
\begin{array}{ll}
Q_1=x^2+y^2+z^2+w^2,&Q_2=xw+yz,\\
Q_3=xz+yw,&Q_4=x^2+y^2-z^2-w^2,\\
Q_5=x^2-y^2-z^2+w^2,&Q_6=x^2-y^2+z^2-w^2,\\
Q_7=xw-yz,&Q_8=xy+zw,\\
Q_9=xy-zw,&Q_{10}=xz-yw.
\end{array}
$$
We use $Q_j$ also for its zero locus when this causes no confusion. These quadrics are smooth. Each contains twelve Klein lines, six in each ruling, and their thirty-six intersection points are $Q_j\cap Z_{31}$.
The intersections of any two distinct fundamental quadrics are unions of four Klein lines. Each Klein line lies on four fundamental quadrics, and each Klein point lies on six.

The incidence counts, including the subdivision induced by $K$, are
\begin{center}
\begin{tabular}{lrrr}
\hline
 & $Z_{31}$ & $Z_{29}$ & $R_{20}$\\
\hline
Points on a Klein line &6&4&2\\
Points on a fundamental quadric &36&24&12\\
\hline
\end{tabular}
\end{center}
The last two columns follow directly from \eqref{eq:z29-coordinates}.
Thus $Q_j\cap Z_{31}$ is a $(6,6)$ grid, and $Q_j\cap Z_{29}$ consists of twenty-four points, four on each of the twelve ruling lines. 
The quadric  $Q_7$ has the standard Segre parametrization:
$$
 ([s_0:s_1],[r_0:r_1])\longmapsto
 [s_0r_0:s_0r_1:s_1r_0:s_1r_1].
$$
\begin{remark}\label{rem:quadrics union and multiplicity}
For later multiplicity arguments, we observe that the product
$\Phi_{20}=\prod_{j=1}^{10}Q_j$ 
has degree $20$ and order $6$ at every point of $Z_{31}$.    
\end{remark}

\subsection{Heisenberg quartics, Maschke coordinates, and the octic}

Let $H_{2,2}\subset\SL_4(\C)$ be generated by
$$
\begin{array}{ll}
(x,y,z,w)\mapsto(z,w,x,y),&(x,y,z,w)\mapsto(y,x,w,z),\\
(x,y,z,w)\mapsto(x,y,-z,-w),&(x,y,z,w)\mapsto(x,-y,z,-w).
\end{array}
$$
It has order $32$, with centre $\{\pm I_4\}$. 
Its projective image
$\cH\simeq(\Z/2\Z)^4$ 
is normal in both projective groups.
Each nonidentity element of $\cH$ fixes two skew lines, and these fifteen pairs are $\call$.

The space $W=\C[x,y,z,w]_4^{H_{2,2}}$ of invariant quartics has as a basis
\begin{equation}\label{eq:heisenberg-basis}
g_0=x^4+y^4+z^4+w^4,\;
g_1=2(x^2y^2+z^2w^2),\;
g_2=2(x^2z^2+y^2w^2),\;
g_3=2(x^2w^2+y^2z^2),\; 
g_4=4xyzw.
\end{equation}
Following the conventions of \cite[Section 2]{Eklund2010}, we consider
\begin{equation}\label{eq:maschke-coordinates}
\begin{array}{lll}
 t_0=\frac13g_0-g_1-g_2-g_3,&
 t_1=\frac13g_0-g_1+g_2+g_3,&
 t_2=\frac13g_0+g_1-g_2+g_3,\\
 t_3=\frac13g_0+g_1+g_2-g_3,&
 t_4=-\frac23g_0+2g_4,&
 t_5=-\frac23g_0-2g_4.
\end{array}
\end{equation}
They span $W$ and satisfy $\sum t_j=0$. In this linear model, $G$ permutes
the six polynomials $t_j$ and the resulting quotient is
$$
 \overline G/\cH\simeq \mathfrak S_6.
$$
The four generators of $K$ in \eqref{eq:generators-of-K} induce, respectively, the transpositions $(2\ 3)$, $(1\ 2)$, $(0\ 1)$ and $(0\ 5)$ on the six indices, while $r_{25}$ induces $(4\ 5)$.
Thus
$$
\overline K/\cH\simeq\mathfrak S_5,
\;\; a_4=-\frac32t_4,
$$
and $K$ is the stabilizer of $t_4$. 
The index-six relation is the stabilizer relation for one of the six Maschke quartics. The six conjugate subgroups correspond to the six choices of distinguished index.

Let $p_r=\sum_{j=0}^5t_j^r$. 
Direct expansion gives the classical Maschke octic, the unique (up to scalar) degree eight $G$-invariant form: 
\begin{equation}\label{eq:maschke-octic}
\begin{aligned}   
f_8&=
x^8+y^8+z^8+w^8+14(x^4y^4+x^4z^4+x^4w^4+y^4z^4+y^4w^4+z^4w^4)+168x^2y^2z^2w^2\\
&=g_0^2+3(g_1^2+g_2^2+g_3^2)+6g_4^2
=\frac34p_2.
\end{aligned}
\end{equation}
Since $G=G_{31}$ has basic invariant degrees $8,12,20,24$, and hence no invariants of degree $4$, the space of degree twelve $G$-invariant forms is one-dimensional; we normalize its generator as 
$f_{12}=-\frac94p_3$.
Explicitly, in the $x,y,z,w$ coordinates it is
\begin{equation}\label{eq:f12-normalization}
\begin{aligned}
f_{12}={}&x^{12}+y^{12}+z^{12}+w^{12}\\
&-33\bigl(
x^8y^4+x^8z^4+x^8w^4
+x^4y^8+x^4z^8+x^4w^8\\
&\hspace{3.5em}
+y^8z^4+y^8w^4+y^4z^8+y^4w^8
+z^8w^4+z^4w^8\bigr)\\
&+330\bigl(
x^4y^4z^4+x^4y^4w^4
+x^4z^4w^4+y^4z^4w^4\bigr)\\
&+792x^2y^2z^2w^2(x^4+y^4+z^4+w^4).
\end{aligned}
\end{equation}

\subsection{The quotient map to the Igusa quartic}

The basis \eqref{eq:heisenberg-basis} has no common zero in $\P^3$.
Indeed, $g_4=0$ forces a coordinate to vanish, after which $g_1=g_2=g_3=0$
forces at most one coordinate to be non-zero; $g_0$ then cannot vanish.
It follows that the forms $t_0,\ldots,t_5$ define a morphism
$$
 \mu:\P^3\longrightarrow\P(W^*),\;\mbox{ with }\;
 P\longmapsto[t_0(P):\cdots:t_5(P)].
$$
We write $T_0,\ldots,T_5$ for the coordinates dual to $t_0,\ldots,t_5$, which identifies $\P(W^*)$ with the hyperplane $\sum T_j=0$ in $\P^5$.
The image of $\mu$ is the Igusa quartic
\begin{equation}\label{eq:igusa}
\cali=\left\{\sum T_j=0,\;\;
\left(\sum T_j^2\right)^2-4\sum T_j^4=0\right
\} \subset\P^5.
\end{equation}
In particular, substitution in the second equation gives the polynomial identity $p_2^2=4p_4$. 
The morphism realizes the quotient by $\cH$ and has degree $16$ \cite[Proposition 2.3]{Eklund2010}. Its degree also agrees with $4^3=16\deg\cali$. 
Since 
$\mu^*\O_{\cali}(1)=\O_{\P^3}(4)$
is ample, no curve is contracted and the morphism is finite.

The singular lines of $\cali$ are indexed by the fifteen partitions of $\{0,\ldots,5\}$ into three pairs. 
For a partition
$\{a,b\}\sqcup\{c,d\}\sqcup\{e,f\}$, the corresponding line is
$$
T_a=T_b,\; T_c=T_d,\; T_e=T_f,\; \sum T_j=0.
$$
The two fixed lines of any non-identity element of $\cH$ have the same image, one of these singular lines. 
The fifteen triple-intersection points are the permutations of 
$[1:1:1:1:-2:-2]$. 
Let $q_{ab}$ denote the point whose two entries $-2$ occupy positions $a,b$. 
The coordinate formulas give
\begin{equation}\label{eq:images-of-points}
\begin{aligned}
 \mu(Z_{31})&=\{q_{ab}:0\le a<b\le5\},\\
 \mu(Z_{29})&=\{q_{ab}:0\le a<b\le5,\ 4\notin\{a,b\}\},\\
 \mu(R_{20})&=\{q_{ab}:0\le a<b\le5,\ 4\in\{a,b\}\}.
\end{aligned}
\end{equation}
For each $a<b$, the fibre $\mu^{-1}(q_{ab})$ consists of four points of $Z_{31}$.
This identifies the subgroup decomposition with the decomposition of the fifteen unordered pairs into ten avoiding the distinguished index $4$ and five containing it.
Every singular line contains three of the $q_{ab}$, one for each pair in its partition. Exactly two avoid the distinguished index, explaining the four-versus-two distribution on each Klein line.

\subsection{The dual Segre cubic and the invariant dodecics}

The space $\P(W)$ parametrizes the Heisenberg invariant quartic surfaces and is dual to the projective space containing $\cali$.
Indeed, since the quartics $t_0,\ldots,t_5$ span $W$ and satisfy the unique linear relation $\sum t_j=0$, every element of $W$ can be written as $\sum u_jt_j$. 
Adding the same constant to all six coefficients does not change this polynomial. 
We remove this ambiguity by requiring $\sum u_j=0$. Then the resulting coefficients are unique, up to a common non-zero scalar when we pass to $\P(W)$.

A parameter $u=[u_0:\cdots:u_5]$ therefore determines both the hyperplane
$$
H_u=\left\{\sum_{j=0}^5u_jT_j=0\right\}
$$
in the target of $\mu$ and the quartic surface
$$
X_u=V\left(\sum_{j=0}^5u_jt_j\right)\subset\P^3.
$$
These constructions are related by
$X_u=\mu^{-1}(H_u\cap\cali)$.
Thus the geometry of the family of invariant quartics is encoded by the hyperplane sections of the Igusa quartic, as expected.

The projective dual $\cali^\vee$ is by definition the closure of the set of hyperplanes tangent to $\cali$ at its smooth points.
In our coordinates it is the Segre cubic
\begin{equation}\label{eq:segre}
\mathcal{S}=
\left\{\sum_{j=0}^5u_j=0,\;\; \sum_{j=0}^5u_j^3=0\right\} \subset\P(W).
\end{equation}
For a general $u\in\mathcal{S}$, the corresponding surface $X_u$ is a Kummer quartic. It has sixteen ordinary double points, forming one orbit under the projective Heisenberg group.
Geometrically, the tangency point of $H_u$ with $\cali$ lifts under $\mu$ to these sixteen points.
Indeed, over a general smooth point of $\cali$, the quotient map $\mu$ is locally an isomorphism at each of the sixteen points in its fiber, so the singularity of the tangent hyperplane section lifts to each of them.

The Segre cubic has ten nodes, represented by the permutations
of 
$$
[1:1:1:-1:-1:-1].
$$
A sign pattern and its simultaneous negative represent the same projective point, giving $\binom{6}{3}/2=10$ nodes.
These special parameters correspond to the ten double quadrics $V(Q_j^2)$, that is, the fundamental quadrics taken with multiplicity two. For example,
$$
Q_7^2=\frac{g_3-g_4}{2},
\;\;
-t_0+t_1+t_2-t_3-t_4+t_5=4(g_3-g_4)=8Q_7^2.
$$
Consequently, $V(Q_7^2)$ corresponds to the node
$[-1:1:1:-1:-1:1]$.

Finally, the fifteen planes contained in $\mathcal{S}$ are indexed by the partitions of $\{0,\ldots,5\}$ into three unordered pairs. 
For example, the partition
$\{\{a,b\},\{c,d\},\{e,f\}\}$ 
gives the plane
$$
u_a+u_b=u_c+u_d=u_e+u_f=0.
$$
These equations imply both equations defining $\mathcal{S}$, since the coordinates cancel in pairs, and so do their cubes.
There are $6!/(2^3\,3!)=15$ such partitions. The indexing set is the same as that labelling the singular lines of $\cali$ above.
For the classical geometry of this family and its duality, see again \cite[Section 2.4]{Eklund2010}.

There is a connection with the invariant dodecics provided by the Gauss map. Let $P_r(T)=\sum T_j^r$. 
Taking the gradient of the quartic equation in \eqref{eq:igusa} and projecting to the hyperplane $\sum u_j=0$ gives
\begin{equation}\label{eq:gauss}
\gamma:\cali\rationalto\mathcal{S},\;\;
T\longmapsto[U_0(T):\cdots:U_5(T)],\;\;
U_j(T)=6T_j^3-\frac32P_2(T)T_j-P_3(T).
\end{equation}
Indeed, the ambient gradient is 
$4P_2T_j-16T_j^3$. Subtracting its mean and multiplying by $-3/8$ gives the formula displayed in \eqref{eq:gauss}. 
The coordinate sum of the $U_j$ is zero.
Their common zero set on $\cali$ is $\operatorname{Sing}\cali$, so the Gauss map is defined on the smooth locus of $\cali$.

We define the six degree twelve forms as follows:
\begin{equation}\label{eq:dodecics}
B_j=\frac1{48}U_j(t)
=\frac18t_j^3-\frac1{32}p_2t_j-\frac1{48}p_3,
\;\text{ for }\; 0\le j\le5.
\end{equation}
Thus $[B_0:\cdots:B_5]$ represents $\gamma\circ\mu$, and the
forms are permuted by $G$. Since $K$ stabilizes $t_4$, the
form $B_4$ is $K$-invariant.
Using $t_4=-\frac23a_4$, $p_2=\frac43f_8$, and
$p_3=-\frac49f_{12}$, we obtain
\begin{equation}\label{eq:b12-formula}
B_4=\frac{-4a_4^3+3a_4f_8+f_{12}}{108}.
\end{equation}
The two Segre relations become 
$\sum B_j=\sum B_j^3=0$.
Expansion in the original coordinates, with each sum taken over the distinct monomials obtained by permuting the variables, gives
\begin{equation}\label{eq:bonnafe-b12}
B_0=
\Sigma(x^8y^2z^2)-\Sigma(x^6y^4z^2) +2\Sigma(x^6y^2z^2w^2)-2\Sigma(x^4y^4z^4)+2\Sigma(x^4y^4z^2w^2).
\end{equation}
Set $B_{12}=B_4$ and $S_{12}=V(B_{12})$, as in the introduction.
The polynomial $B_0$ in \eqref{eq:bonnafe-b12} equals to the polynomial $f_3$ displayed in \cite[Section~1.B]{Bonnafe2025Dodecic}, with its last variable $t$ renamed $w$.
To compare it with $B_{12}$, consider the involution
$$
C=r_{30}=\frac12
\begin{pmatrix}
 1&-1&-1&1\\
 -1&1&-1&1\\
 -1&-1&1&1\\
 1&1&1&1
\end{pmatrix}\in G.
$$
Writing $X=(x,y,z,w)^T$, direct substitution gives $t_0(CX)=t_4(X)$, so $CKC^{-1}$ is the stabilizer of $t_0$ in $G$.
Since $p_2$ and $p_3$ are $G$-invariant, \eqref{eq:dodecics} then gives
\begin{equation}\label{eq:bonnafe-comparison}
B_{12}(X)=B_4(X)=B_0(CX).
\end{equation}
Thus the projective transformation induced by $C$ identifies $S_{12}$ with the surface $V(B_0)$ studied by Bonnaf\'e.
In the two-parameter family of dodecics considered in \cite[Section~1.B]{Bonnafe2025Dodecic}, $V(B_0)$ is the member with both parameters equal to $0$.
It is irreducible, and its singular locus is the union of the thirty Klein lines \cite[Lemma~1.2 and Section~2.C(f)]{Bonnafe2025Dodecic}.
Since $C$ permutes these lines, the same conclusions hold for $S_{12}$.
The surface had already appeared in \cite[Example~5.5(4)]{Bonnafe2020}.

For a non-zero homogeneous form $F$ vanishing at a point $p$, we denote by $\Loc_p(F)$ the lowest degree non-zero homogeneous part of its dehomogenization in the chosen affine coordinates centred at $p$.
This polynomial defines the tangent cone of $V(F)$ at $p$, and its degree equals $\mult_p(F)$.
\begin{proposition}\label{prop:b12-orders}
The form $B_{12}$ has order four at every point of $Z_{29}$
and order three at every point of $R_{20}$.
Up to scalar, it is the unique $K$-invariant form of degree
twelve vanishing to order at least four at $Z_{29}$.
\end{proposition}

\begin{proof}
Since $B_{12}$ is $K$-invariant and $Z_{29},R_{20}$ are
$K$-orbits, it suffices to inspect one point from each orbit.

At $P_{25}=[1:0:0:0]\in R_{20}$, in the chart $x=1$,
the lowest-degree term is
$$
\Loc_{P_{25}}(B_{12})=yzw.
$$
Hence the order at $P_{25}$ is three.

At $P_{22}=[1:i:0:0]\in Z_{29}$, putting $y=i+v$,
the lowest-degree part is
$$
\Loc_{P_{22}}(B_{12})
=4(z^2+w^2)^2-16iv^2zw.
$$
This non-zero quartic shows that the order at $P_{22}$
is four.

For the uniqueness, the basic invariant degrees $4,8,12,20$
give
$\dim_{\C}\C[x,y,z,w]_{12}^{K}=3$.

On the line $z=w=0$, expansion at $P_{22}$ gives
$$
\begin{aligned}
a_4f_8(1,i+v,0,0)
  &=32-192iv-576v^2+O(v^3),\\
a_4^3(1,i+v,0,0)
  &=8-48iv-168v^2+O(v^3).
\end{aligned}
$$
Together with the order-four vanishing of $B_{12}$,
these expansions show that
$B_{12},a_4f_8,a_4^3$ are linearly independent.
They therefore form a basis of the degree-twelve
invariant space.

If $F=c_0B_{12}+c_1a_4f_8+c_2a_4^3$ vanishes to order at least four at $P_{22}$,
its restriction to $z=w=0$ has zero constant and
quadratic coefficients. Thus
$$
32c_1+8c_2=0,
\;\mbox{ and }\;
576c_1+168c_2=0.
$$
Since $32\cdot168-8\cdot576=768\neq0$, these equations force $c_1=c_2=0$, proving uniqueness.
\end{proof}
\begin{proposition}[Dodecics from Igusa--Segre duality]
\label{prop:igusa-segre-dodecics}
The rational map $\gamma\circ\mu:\P^3\dashrightarrow\mathcal S$ is represented by $[B_0:\cdots:B_5]$.
For each $j$, the dodecic $V(B_j)$ is the closure in $\P^3$ of the pullback of the coordinate hyperplane section $\mathcal S\cap V(u_j)$ on the domain of definition of this map.
These six dodecics form the $G$-orbit of $S_{12}=V(B_4)$.
Their defining forms have the $G$-invariant orbit product of degree $72$
\begin{equation}\label{eq:degree72-norm}
B_{72}=\prod_{j=0}^5B_j,\;\;
\mult_P(B_{72})=22\;\;
\text{for every }P\in Z_{31}.
\end{equation}
In particular, $S_{72}=V(B_{72})$ gives $\alphahat(Z_{31})\leq72/22=36/11$.
\end{proposition}

\begin{proof}
The identity $U_j(t)=48B_j$ in \eqref{eq:dodecics} gives the asserted coordinates of $\gamma\circ\mu$.
Their common zero set is $\mu^{-1}(\operatorname{Sing}\cali)$, which has codimension two because $\mu$ is finite and the singular locus of $\cali$ is a union of lines.
Consequently, taking the closure of each pullback gives the divisor $V(B_j)$.
The action of $G$ permutes the six indices transitively, so the six forms are the translates of $B_4$ and their product is $G$-invariant.

Let $P\in Z_{31}$ and write $\mu(P)=q_{ab}$.
By \eqref{eq:images-of-points}, equivariance and Proposition~\ref{prop:b12-orders}, the order of $B_j$ at $P$ is four when $j\notin\{a,b\}$ and three when $j\in\{a,b\}$.
There are four indices of the first kind and two of the second, so the order of the product at $P$ is $4\cdot4+2\cdot3=22$.
Its degree is $6\cdot12=72$, giving the stated upper bound.
\end{proof}

\subsection{The Waldschmidt constant of the Klein configuration}
The lower bound on $\alphahat(Z_{31})$ follows from
the incidence geometry of the Klein lines and fundamental quadrics.

\begin{lemma}\label{lem:Klein lower bound}
We have
$$
\alphahat(Z_{31})\geq\frac{36}{11}.
$$
\end{lemma}
\begin{proof}
Suppose, to the contrary, that there is a surface $S$ of
degree $d$ vanishing to order at least $m$ at every point
of $Z_{31}$, where $m$ is a positive integer and
$11d<36m$. Among all such pairs $(S,m)$, choose one with $d$ minimal.

Let
$a=\min_{\ell}\mult_\ell(S)$,
where the minimum is taken over the $30$ Klein lines,
and choose a line $\ell$ attaining this minimum.
For a general plane $\Pi$ through $\ell$, the line $\ell$
occurs in the trace $S|_\Pi$ with multiplicity equal $a$.
Thus the residual curve
$$
C=S|_\Pi-a\ell
$$
has degree $d-a$, does not contain $\ell$, and has
multiplicity at least $m-a$ at each of the six points
of $Z_{31}$ on $\ell$. By B\'ezout,
\begin{equation}\label{eq:lower bound on a}
d-a\geq 6(m-a),
\;\mbox{ hence }\;
5a\geq 6m-d.
\end{equation}

Suppose that a fundamental quadric $Q$ is not a component
of $S$. The trace $\Gamma=S|_Q$ has bidegree $(d,d)$.
The twelve Klein lines on $Q$ form a $6\times6$ grid,
and each occurs in $\Gamma$ with multiplicity at least $a$.
Subtracting $a$ times their union leaves an effective
divisor of bidegree $(d-6a,d-6a)$. Consequently,
$$
d\geq 6a\geq 6\frac{6m-d}{5},
$$
which gives $11d\geq36m$, a contradiction.

Therefore all ten fundamental quadrics are components
of $S$. In particular, $d\geq20$, and $11d<36m$
then implies $m>6$.
By Remark~\ref{rem:quadrics union and multiplicity},
the union of all $10$ fundamental quadrics has degree $20$ and multiplicity $6$
at every point of $Z_{31}$.
Subtracting this union from $S$ gives a residual surface
$S'$ of degree
$d'=d-20$
and multiplicity at least $m'=m-6>0$ at every point of $Z_{31}$. Moreover,
$$
36m'-11d'=36(m-6)-11(d-20)=(36m-11d)+4>0.
$$
Thus $S'$ satisfies $11d'<36m'$. Since its degree is smaller than that of $S$, this
contradicts the choice of $S$.
\end{proof}
Combining Lemma \ref{lem:Klein lower bound} with Proposition \ref{prop:igusa-segre-dodecics} proves Theorem~B.

\subsection{The Waldschmidt constant of \texorpdfstring{$Z_{29}$}{Z29}}
\label{subsec:G29-Waldschmidt}

The dodecic $B_{12}$ gives an upper bound of $3$ for $\alphahat(Z_{29})$.
For the matching lower bound, we restrict to the twenty Klein planes corresponding to $R_{20}$.
Let $\mathcal H_{20}=\{H_P:P\in R_{20}\}$, with the Hermitian polarity fixed in the introduction.
Since $R_{20}$ is a single $K$-orbit and $K$ acts unitarily, the identity $H_{gP}=g(H_P)$ shows that $\mathcal H_{20}$ is also a single $K$-orbit.
It contains the coordinate plane $H_0=H_{P_{25}}=V(x)$.

\begin{lemma}\label{lem:G29-plane-sections}
Let $H\in\mathcal H_{20}$ and let $C\subset H$ be a plane curve of degree $e$ vanishing to order at least $n$ at every point of $Z_{29}\cap H$, where $n$ is a positive integer.
Then $e\geq3n$.
\end{lemma}

\begin{proof}
Since $K$ preserves $Z_{29}$ and acts transitively on $\mathcal H_{20}$, it suffices to consider $H_0=V(x)$.
The twelve points of $Z_{29}\cap H_0$ are $P_1,\ldots,P_{12}$.
They lie four at a time on the three Klein lines $V(x,y)$, $V(x,z)$ and $V(x,w)$.
These lines form a triangle whose vertices are $P_{26},P_{27},P_{28}\in R_{20}$, so none of the twelve points is a vertex.

Suppose that $e<3n$.
Each side of the triangle contains four points of multiplicity at least $n$, so B\'ezout forces all three sides to be components of $C$.
Their union has degree three and multiplicity one at each of the twelve points.
Subtracting it leaves an effective curve of degree $e-3$ and multiplicity at least $n-1$ at those points.
The strict inequality between the degree and three times the prescribed multiplicity persists after subtraction.
Repeating this argument $n$ times would leave an effective curve of degree $e-3n<0$, which is impossible.
\end{proof}
Every plane in $\mathcal H_{20}$ contains twelve points of $Z_{29}$.
Since $K$ is transitive on $Z_{29}$, each point belongs to the same number of these planes.
Counting incidences gives $20\cdot12/40=6$ planes through each point.
Consequently, their reduced union $D_{20}$ has degree twenty and multiplicity six at every point of $Z_{29}$.

\begin{proposition}\label{prop:G29-Waldschmidt}
We have
$$
\alphahat(Z_{29})=3.
$$
\end{proposition}

\begin{proof}
The upper bound follows from the degree and order-four vanishing of $B_{12}$ established above.
For the lower bound, suppose that there is a surface $S$ of degree $d$ vanishing to order at least $m$ at every point of $Z_{29}$, where $m$ is a positive integer and $d<3m$.
Choose such a surface with $d$ minimal.

If a plane $H\in\mathcal H_{20}$ were not a component of $S$, the trace $S|_H$ would be a plane curve of degree $d$ vanishing to order at least $m$ at every point of $Z_{29}\cap H$.
Lemma~\ref{lem:G29-plane-sections} would give $d\geq3m$, a contradiction.
Thus all twenty planes are components of $S$, and we may subtract $D_{20}$.

Since $d\geq20$ and $d<3m$, we have $m>6$.
The residual surface $S'=S-D_{20}$ has degree $d'=d-20$ and multiplicity at least $m'=m-6>0$ at every point of $Z_{29}$.
Moreover, $3m'-d'=(3m-d)+2>0$.
Thus $S'$ is a counterexample of smaller degree, contradicting the choice of $S$.
\end{proof}

This proves Theorem~A and completes the calculation for $Z_{29}$ and $Z_{31}$.

\section{The Penrose configuration}
\label{sec:penrose}

We now consider the forty-point configuration $Z_{32}$ associated with $G_{32}$.
The lower bound again comes from restriction to distinguished planes.
On each such plane, two related Hesse arrangements yield a weighted B\'ezout estimate, while the four-point lines impose additional vanishing at nine auxiliary points.
Together, these conditions force every polar plane to be a component of a hypothetical surface whose degree-to-multiplicity ratio is less than $10/3$, and subtraction of their union gives a descent analogous to those used in the preceding section.

\subsection{The configuration and its polar planes}

We use homogeneous coordinates $[x:y:z:w]$ on $\P^3$ and the realization of the Penrose configuration described in \cite{ChiantiniEtAl2022}.
Its homogeneous ideal is
\begin{equation}\label{eq:penrose-ideal}
\begin{aligned}
I(Z_{32})=\bigl(&xyzw,\ w(x^3-y^3+z^3),\ z(x^3+y^3+w^3),\\
               &y(-x^3+z^3+w^3),\ x(y^3+z^3-w^3)\bigr).
\end{aligned}
\end{equation}
Let $\omega$ be a primitive cube root of unity, so that $\overline\omega=\omega^2$.
The configuration consists of the four coordinate points and the following four families of nine points:
\begin{equation}\label{eq:penrose-coordinates}
\begin{array}{ll}
[0:1:\omega^a:-\omega^b], & [1:0:-\omega^a:-\omega^b],\\[2pt]
[1:\omega^a:0:\omega^b], & [1:-\omega^a:\omega^b:0],
\end{array}
\;\; a,b\in\{0,1,2\}.
\end{equation}
This configuration goes back to Zimba and Penrose \cite{ZimbaPenrose1993}; it is arithmetically Gorenstein and $(5,8)$-geproci \cite{ChiantiniEtAl2022}.

For $A=[a:b:c:d]\in Z_{32}$, its Hermitian polar plane is $H_A=A^\perp=V(\overline a x+\overline b y+\overline c z+\overline d w)$.
In our unitary realization of $G_{32}$, the points of $Z_{32}$ are the normal directions to the forty reflecting hyperplanes, so the planes $H_A$ are precisely their projectivizations.
The configuration is stable under complex conjugation, as is clear from its displayed coordinates, and the ordinary dual configuration is therefore represented by the same coordinate list.
Hermitian polarity is equivariant: $H_{gA}=g(H_A)$ for $g\in G_{32}$.
We use the description of the reflection arrangement and its stabilizers given in \cite[Section~5.4]{CousinMoussard2018}.

\begin{proposition}\label{prop:penrose-global-incidences}
Every plane $H_A$ contains twelve points of $Z_{32}$.
Every point of $Z_{32}$ belongs to twelve of these planes.
There are ninety lines containing four points of $Z_{32}$, no line contains five or more, and every point belongs to nine four-point lines.
The group $G_{32}$ acts transitively on the forty points, the forty planes and the ninety four-point lines.
\end{proposition}
\begin{proof}
By \cite[Proposition~5.4.7]{CousinMoussard2018}, the rank-two flats of the reflection arrangement form two orbits: $240$ intersections of two orthogonal hyperplanes and $90$ intersections of four nonorthogonal hyperplanes.
The group is also transitive on its forty reflecting hyperplanes and hence, by polarity, on $Z_{32}$.
The orthogonal pairs give $2\cdot240/40=12$ points in each polar plane.
The relation $A\in H_B$ is symmetric in $A$ and $B$, so each point also belongs to twelve polar planes.

A line through $r$ configuration points is polar to a projective line contained in the corresponding $r$ reflecting planes.
Conversely, every rank-two flat polarizes to a line through the corresponding normal directions.
The description of the flats therefore gives ninety four-point lines and excludes lines with five or more configuration points.
Transitivity and incidence counting give $90\cdot4/40=9$ four-point lines through each point.
\end{proof}

\subsection{Hesse arrangements on  polar planes}

Fix $A=[v]\in Z_{32}$ and let $H=H_A=\P(v^\perp)$ and $X=Z_{32}\cap H$.
Since the Hermitian form is positive definite, $A\notin H$.
By Steinberg's theorem, the vector stabilizer $(G_{32})_v=\{g\in G_{32}:gv=v\}$ is generated by the reflections fixing $v$ \cite{Steinberg1964}.
Since the action is unitary, these reflections preserve $v^\perp$ and restrict to reflections on this three-dimensional space.
The resulting reflection group is isomorphic to $G_{25}$, and its projective action on $H$ is the Hessian group \cite[Lemma~5.4.4]{CousinMoussard2018}.

Polar points of lines contained in $H$ are taken with respect to the restricted Hermitian form on $v^\perp$.
For a projective line $\ell\subset\mathbb P^3$, the notation $\ell^\perp$ denotes its ambient Hermitian polar line.

Let $\mathcal T$ be the union of the twelve lines $H\cap H_B$ with $B\in X$; these are the reflection lines of $G_{25}$ in $H$.
Let $\Sigma$ be the union of the lines in $H$ that contain four points of $Z_{32}$.
In the terminology of \cite[Section~5.3]{CousinMoussard2018}, the components of $\Sigma$ are the projectivizations of the proper planes of $G_{25}$; this description is not used below.

\begin{proposition}\label{prop:penrose-hesse-incidences}
The set $X$ consists of the twelve double points of $\mathcal T$, and $\Sigma$ consists of nine lines.
The nine four-point lines through $A$ meet $H$ in nine distinct points forming a set $Y$ disjoint from $X$.
These are the nine quadruple points of $\mathcal T$.
Every point of $X$ lies on three components of $\Sigma$ and two components of $\mathcal T$.
Every point of $Y$ lies on four components of $\mathcal T$ and no component of $\Sigma$.
Each component of $\Sigma$ contains four points of $X$.
Each component of $\mathcal T$ contains two points of $X$ and three points of $Y$.
Thus $(X,\Sigma)$ is the dual Hesse configuration $(12_3,9_4)$, while $(Y,\mathcal T)$ is the Hesse configuration $(9_4,12_3)$.
\end{proposition}
\begin{proof}
By Proposition~\ref{prop:penrose-global-incidences}, $G_{32}$ acts transitively on $Z_{32}$, and $H_{gA}=g(H_A)$; hence $g$ carries $X$, $Y$, $\Sigma$ and $\mathcal T$ for $A$ to the corresponding objects for $gA$.
It therefore suffices to take $A=[1:0:0:0]$, so that $H=V(x)$, and we use $[y:z:w]$ as coordinates on $H$.
By \eqref{eq:penrose-coordinates}, $X$ consists of the three coordinate points and the nine points $[1:\omega^a:-\omega^b]$ with $a,b\in\{0,1,2\}$.
For $B\in X$, the line $H\cap H_B$ has equation $\overline B\cdot(y,z,w)^T=0$, so
$$
\mathcal T=V\Bigl(yzw\prod_{a,b=0}^{2}\bigl(y+\omega^az-\omega^bw\bigr)\Bigr).
$$
The line through $A$ and $[1:0:-\omega^a:-\omega^b]$ contains the three points $[1:0:-\omega^{a'}:-\omega^{b'}]$ with $b'-a'\equiv b-a\pmod 3$ and meets $H$ at $[0:1:\omega^{b-a}]$.
Treating the families $[1:\omega^a:0:\omega^b]$ and $[1:-\omega^a:\omega^b:0]$ in the same way, and noting that $A$ lies on no other line through four points of $Z_{32}$, we see that the nine four-point lines through $A$ meet $H$ in the set
$$
Y=\bigl\{[0:1:\omega^c],\ [1:0:\omega^c],\ [1:-\omega^c:0]\;:\;c\in\{0,1,2\}\bigr\}.
$$
Each of the nine lines of
$$
V\bigl((y^3-z^3)(z^3+w^3)(w^3+y^3)\bigr)
$$
passes through one coordinate point and three of the points $[1:\omega^a:-\omega^b]$.
Conversely, a line in $H$ contains at most two coordinate points, and if it contains two it is a coordinate line, which meets none of the nine points $[1:\omega^a:-\omega^b]$.
These nine points lie on the smooth cubic $V(y^3-2z^3-w^3)$, so a line contains at most three of them.
Hence a line in $H$ containing four points of $Z_{32}$ passes through exactly one coordinate point and three of the nine points, and is one of the nine lines above; this identifies $\Sigma$.

All asserted incidences can now be read off from these equations.
For instance, $[1:\omega^a:-\omega^b]$ lies on the line $y+\omega^{a'}z-\omega^{b'}w=0$ if and only if $1+\omega^{a+a'}+\omega^{b+b'}=0$, which happens for exactly two pairs $(a',b')$, and it lies on no coordinate line.
Similarly, $[0:1:\omega^c]$ lies on $V(y)$ and on the three lines $y+\omega^az-\omega^bw=0$ with $a\equiv b+c\pmod 3$, and it does not lie on $\Sigma$.
The remaining incidences are checked in the same way.
\end{proof}

\subsection{A weighted B\'ezout estimate}

The two arrangements give a weighted B\'ezout estimate for plane divisors, which we will apply to restrictions of surfaces to polar planes.

\begin{lemma}\label{lem:penrose-weighted-bezout}
Let $X,Y,\Sigma,\mathcal{T}$ be the point sets and arrangements defined in the preceding proposition.
For every effective plane divisor $C$ with rational coefficients, we have
$$
\deg C \geq
\frac14\sum_{P\in X}\operatorname{mult}_P(C)
+\frac16\sum_{P\in Y}\operatorname{mult}_P(C).
$$
In particular, if $C$ has multiplicity at least $a\geq0$ at every point of $X$ and at least $b\geq0$ at every point of $Y$, then $\deg C\geq3a+\frac32b$.
\end{lemma}

\begin{proof}
Both degree and point multiplicities are additive on effective divisors, so it suffices to prove the estimate for an irreducible curve $D$.
Write $e=\deg D$, $M_X=\sum_{P\in X}\operatorname{mult}_P(D)$ and $M_Y=\sum_{P\in Y}\operatorname{mult}_P(D)$.

Suppose first that $D$ is not a component of either arrangement.
The divisor $\Sigma$ has degree nine and multiplicity three at every point of $X$, while $\mathcal{T}$ has degree twelve and multiplicities two at $X$ and four at $Y$.
B\'ezout's theorem therefore gives
$$
9e\geq3M_X,\;\; 12e\geq2M_X+4M_Y.
$$
Dividing these inequalities by $18$ and $24$, respectively, and adding yields $e\geq\frac14M_X+\frac16M_Y$.

If $D$ is a component of $\Sigma$, then $e=1$, $M_X=4$ and $M_Y=0$, so equality holds.
If $D$ is a component of $\mathcal{T}$, then $e=1$, $M_X=2$ and $M_Y=3$, and again equality holds.
This proves the estimate in all cases.
The final assertion follows from $|X|=12$ and $|Y|=9$.
\end{proof}

\subsection{The Waldschmidt constant}

\begin{theorem}\label{thm:penrose-waldschmidt}
We have
$$
\alphahat(Z_{32})=\frac{10}{3}.
$$
\end{theorem}
\begin{proof}
For each $A\in Z_{32}$, choose a linear equation $h_A$ of its polar plane $H_A$ and put $\Phi_{40}=\prod_{A\in Z_{32}}h_A$.
By Proposition~\ref{prop:penrose-global-incidences}, this form has degree forty and order twelve at every configuration point.
Consequently, $\widehat{\alpha}(Z_{32})\leq40/12=10/3$.

For the lower bound, suppose that an effective surface $S$ of degree $d$ vanishes to order at least a positive integer $m$ at every point of $Z_{32}$ and satisfies $3d<10m$.
Choose such a surface with $d$ minimal.
Lemma~\ref{lem:penrose-forced-line}, applied with $s=4$, shows that the order of $S$ along every four-point line is at least $r=(4m-d)/3>0$.

Fix $A\in Z_{32}$ and let $H=H_A$, with the associated point sets $X,Y$ and arrangements $\Sigma,\mathcal{T}$.
Suppose that $H$ is not a component of $S$.
Then $\Gamma=S|_H$ is an effective plane divisor of degree $d$ and has multiplicity at least $m$ at every point of $X$.
For each $y\in Y$, there is a four-point line $\ell$ through $A$ meeting $H$ at $y$.
Since $A\notin H$, this line is not contained in $H$.
The order of $S$ along $\ell$ is at least $r$, so its restriction to $H$ has multiplicity at least $r$ at $y$.
Applying Lemma~\ref{lem:penrose-weighted-bezout} directly to $\Gamma$ gives
$$
d\geq3m+\frac32r=5m-\frac d2,
$$
which implies $3d\geq10m$, a contradiction.
Thus $H$ is a component of $S$.

Since $A$ was arbitrary, all forty polar planes are components of $S$.
Their reduced union $D_{40}=V(\Phi_{40})$ has degree forty and multiplicity twelve at every point of $Z_{32}$.
The inequalities $d\geq40$ and $3d<10m$ imply $m>12$.
The residual surface $S'=S-D_{40}$ is therefore effective, has degree $d'=d-40$ and has multiplicity at least $m'=m-12>0$ at every configuration point.
Moreover, $10m'-3d'=10m-3d>0$.
Thus $S'$ is a counterexample of smaller degree, contradicting the choice of $S$.
\end{proof}

This proves Theorem~C.
\section{Resurgence and asymptotic resurgence}
\label{sec:resurgence}

For the homogeneous ideal $I$ of a finite set of points, the resurgence
$$
\rho(I)=\sup\left\{\frac mr\;:\;I^{(m)}\not\subseteq I^r\right\}
$$
was introduced in \cite{BocciHarbourne2010}, and the asymptotic resurgence
$$
\widehat\rho(I)=\sup\left\{\frac mr\;:\;I^{(mt)}\not\subseteq I^{rt}\ \text{for all}\ t\gg0\right\}
$$
in \cite{GuardoHarbourneVanTuyl2013}.
Both are controlled by the Waldschmidt constant.
Let $\omega(I)$ denote the largest degree of a minimal homogeneous generator of $I$ and $\operatorname{reg}(I)$ its Castelnuovo--Mumford regularity.
By \cite[Theorem~1.2.1]{BocciHarbourne2010},
\begin{equation}\label{eq:resurgence-bounds}
\frac{\alpha(I)}{\alphahat(I)}\leq\rho(I)\leq\frac{\operatorname{reg}(I)}{\alphahat(I)},
\end{equation}
and by \cite[Lemma~2.3.2(b) and Corollary~2.3.9]{BocciHarbourne2010},
\begin{equation}\label{eq:asymptotic-resurgence-bounds}
\frac{\alpha(I)}{\alphahat(I)}\leq\widehat\rho(I)\leq\frac{\omega(I)}{\alphahat(I)};
\end{equation}
see also \cite{GuardoHarbourneVanTuyl2013} and \cite[(1) and (2)]{Negative_curves_IMRN2019}.
In particular, $\widehat\rho(I)=\alpha(I)/\alphahat(I)$ whenever $I$ is generated in a single degree.

\begin{corollary}\label{cor:resurgence}
Let $I_{29}$, $I_{31}$ and $I_{32}$ be the homogeneous ideals of $Z_{29}$, $Z_{31}$ and $Z_{32}$. Then
$$
\widehat\rho(I_{29})=\frac53,\;\;
\widehat\rho(I_{31})=\frac{11}6,\;\;
\widehat\rho(I_{32})=\frac65,
$$
and
$$
\frac53\leq\rho(I_{29})\leq2,\;\;
\frac{11}6\leq\rho(I_{31})\leq\frac{77}{36},\;\;
\frac65\leq\rho(I_{32})\leq\frac{12}5.
$$
\end{corollary}

\begin{proof}
A direct computation shows that $I_{29}$ is minimally generated by sixteen quintics and $I_{31}$ by twenty-four sextics, and that the $h$-vectors of $Z_{29}$ and $Z_{31}$ are $(1,3,6,10,15,5)$ and $(1,3,6,10,15,21,4)$.
The ideal $I_{32}$ is generated by the five quartics in \eqref{eq:penrose-ideal}, and $Z_{32}$ has $h$-vector $(1,3,6,10,10,6,3,1)$.
Hence $\alpha(I)=\omega(I)$ equals $5$, $6$ and $4$, and $\operatorname{reg}(I)$ equals $6$, $7$ and $8$, for $I=I_{29}$, $I_{31}$ and $I_{32}$, respectively.
The assertions follow from \eqref{eq:resurgence-bounds}, \eqref{eq:asymptotic-resurgence-bounds} and Theorems~A--C.
\end{proof}

For comparison, in rank three the resurgence of the ideal of the singular points of the Klein and of the Wiman line arrangements, associated with $G_{24}$ and $G_{27}$, equals $3/2$, while their asymptotic resurgences lie in $[16/13,816/661]$ and equal $32/27$, respectively \cite[Theorems~1.4 and~1.5]{Negative_curves_IMRN2019}.
\medskip 

\paragraph{Acknowledgements.}
This work grew out of a working group at the Mini-Workshop on \emph{Asymptotic Invariants of Linear Systems and Homogeneous Ideals with Symmetries}, held at the Mathematisches Forschungsinstitut Oberwolfach in September 2026.
We thank the institute for its hospitality and excellent working conditions.
\raggedbottom

\end{document}